\documentclass[12pt]{amsart}

\usepackage{color}
\usepackage{amsfonts, amsmath, amsfonts, amssymb}
\usepackage{graphicx}
\usepackage{algorithm,algpseudocode}
\usepackage{MnSymbol,wasysym}
\definecolor{RedClr}{rgb}{1,0,0}
\definecolor{BlueClr}{rgb}{0,0,1}
\definecolor{TextColor}{rgb}{0,0,0.5}
\definecolor{Violet}{rgb}{0.5,0,1}
\definecolor{Bordeaux}{rgb}{1,0.3,0.4}

\begin{document}

\newtheorem{thm}{Theorem}[section]
\newtheorem{cor}[thm]{Corollary}
\newtheorem{lem}[thm]{Lemma}
\newtheorem{prop}[thm]{Proposition}
\theoremstyle{definition}
\newtheorem{defn}[thm]{Definition}
\theoremstyle{remark}
\newtheorem{rem}[thm]{Remark}
\numberwithin{equation}{section} \theoremstyle{quest}
\newtheorem{quest}[]{Question}
\numberwithin{equation}{section} \theoremstyle{prob}
\newtheorem{prob}[]{Problem}
\numberwithin{equation}{section} \theoremstyle{answer}
\newtheorem{answer}[]{Answer}
\numberwithin{equation}{section}
\theoremstyle{fact}
\newtheorem{fact}[thm]{Fact}
\numberwithin{equation}{section}
\theoremstyle{facts}
\newtheorem{facts}[thm]{Facts}
\numberwithin{equation}{section}
\newtheorem{conj}[]{Conjecture}
\newtheorem{cntxmp}[thm]{Counterexample}
\numberwithin{equation}{section}
\newtheorem{exmp}[thm]{Example}
\numberwithin{equation}{section}
\newtheorem{ex}[thm]{Exercise}
\numberwithin{equation}{section}
\newenvironment{prf}{\noindent{\bf Proof}}{\\ \hspace*{\fill}$\Box$ \par}
\newenvironment{skprf}{\noindent{\\ \bf Sketch of Proof}}{\\ \hspace*{\fill}$\Box$ \par}

\newcommand{\convGH}{\raisebox{-0.5em}{$\stackrel{\longrightarrow}{\scriptstyle GH}$}}

\title[Forman-vs-PH]{Discrete Morse Theory, Persistent Homology and Forman-Ricci Curvature}
\author{Emil Saucan}

\address{Department of Applied Mathematics, ORT Braude College, Karmiel, Isarel}

\email{semil@braude.ac.il}%


\date{\today}


\maketitle


\begin{abstract}
	Using Banchoff's discrete Morse Theory, in tandem with Bloch's result on the strong connection between the former and  Forman's Morse Theory, and our own previous algorithm based on the later, we show that there exists a curvature-based, efficient Persistent Homology scheme for networks and hypernetworks. We also broaden the proposed method to include more general types of networks, by using Bloch's extension of Banchoff's work.  
	Moreover, we show the connection between defect and Forman's Ricci curvature that exists in the combinatorial setting, thus explaining previous empirical results showing very strong correlation between Persistent Homology results obtained using Forman's Morse Theory on the one hand, and Forman's Ricci curvature, on the other.
\end{abstract}


\section{Introduction}

The present paper  is motivated by a number of the author's converging research interests. The basic stimulus is, to be sure, the ever growing importance in a variety of applied fields of the Persistent Homology \cite{ELZ}, \cite{Ca}, \cite{CZ}, \cite{CISZ}, \cite{SEH}, not the least among them being the intelligence of Complex Networks \cite{HMR}, \cite{Petri+}, \cite{PSDV} -- to name just a few among many others; combined with our own sustained interest in Forman's discrete Ricci curvature \cite{Fo}. Of special importance here is our observation \cite{SJ} that, in certain aspects and instances Forman's Ricci curvature might complement (and perhaps even supplement, for some specific tasks), the Persistent Homology method. Another incentive has come from our own recent application \cite{KSRS} of a well known discrete Morse Theory, again due to Forman \cite{Fo98}, \cite{Fo98a}, to the study of Complex Networks, and the further empirical observation \cite{RVRS} that the results provided by the topological approach correlate almost up to coincidence with those produced using a Forman-Ricci curvature scheme.  
These topological and geometric ideas coalesce naturally to our own observation \cite{SJ}, \cite{WSJ1} that a better understanding of networks can only gain by viewing them -- and multiplex networks a fortiori -- as polyhedral complexes, hence all geometric characteristics of such geometric objects, and in particular curvature, can be applied to their study.  (In fact, a trend of modeling networks as simplicial complexes is prevailing these days -- see, e.g. \cite{CB}, \cite{CB1},\cite{DeSG}, ֿ\cite{IPBL}, \cite{K-GB}, \cite{KSRS}, \cite{PSDV}, \cite{PSDV}, \cite{R++}, \cite{SGB}, \cite{WMRB}.) Furthermore, we have show in \cite{SW19} that hypernetworks can be naturally modeled as polyhedral complexes, endowed with natural notions of curvature. Moreover, we recently made the observation that that there exists a canonical way of viewing hyperneworks as simplicial complexes \cite{SS}, hence once again inherent notions of curvature are applicable. Moreover, curvature is strongly related to (and in fact it defines) a polyhedral Morse Theory due to Banchoff \cite{Ban67}, \cite{Ban70}, \cite{Ban83}. Given the surprising fact, proven by Bloch \cite{Bloch04}, that the widely divergent in their definition (and setting) Morse Theories of Forman and Banchoff, are essentially interchangeable, it is only naturally to conclude that the somewhat counterintuitive Forman-Morse approach to Persistent Homology can be replaced, at least in low dimensions, by a simpler curvature-driven one. This method can be extended to a larger, and even better fitting for the modeling of hypernetworks, class of simplicial complexes, using Bloch's generalization of Bachoff's work \cite{Bloch04}. Moreover, we show that, for combinatorial simplicial and more general polyhedral complexes, there is a strong connection between Forman's Ricci curvature and the defect definition of curvature employed in the works of Banchoff and Bloch, a fact which explains the above mentioned correlation observed in \cite{RVRS}.

The reminder of the paper is structured in a natural manner as follows: In Section 2 we present Banchoff's Morse Theory,  followed in Section 3 by an overview of Bloch's proof of the Forman-Banchoff Morse Theories quasi-equivalence. Section 4 is dedicated to the  introduction of Bloch's generalization of Banchoff's Morse Theory. In Section 5 we show the strong connection between angle defect and Forman-Ricci curvature. We conclude in Section 6 with a review of the ideas and results expounded in the previous sections and an overview of the future tasks and list of problems we deem more immediate.

While on the one hand we certainly would wish to ensure that the present paper is self-contained, on the other hand we  would also like the paper not to expand inordinately. To try and satisfy both of these contradictory demands we do not dwell into any details of the Forman-Morse Theory and its application to the Persistent Homology of networks, since this is a subject familiar, we believe, to many of the potential readers and, furthermore, we have quite recently expounded on both the theoretical background and it desired application in our article \cite{KSRS}. Furthermore, we do not present any details regarding Forman's Ricci curvature, beyond the strictly necessary ones, because we did explicate on this notion, especially in its utility in the study of networks and hypernetworks, in a number of papers \cite{SMJSS}, \cite{WSJ1}, \cite{WSJ}, \cite{SSGLSJ}, \cite{SSWJ}, \cite{SW19}. Also, not to repeat ourselves, we do not review here the modeling of hypernetworks as polyhedral complexes devised in \cite{SW19}, nor the canonical view of these structures as simplicial complexes that is shortly forthcoming in \cite{SS}.
\subsection*{Acknowledgment}

The author is indebted to Indrava Roy and Areejit Samal for rekindling his interest in discrete Morse Theory. 


\section{Banchoff's Polyhedral Morse Theory}

This section is dedicated to a review of Banchoff's polyhedral Morse Theory. We first present the more intuitive case of polyhedral surfaces which we expose in a simple, detailed manner, since this will permit, in our opinion, a better future implementation by any interested reader. We follow it by the exposition of the generalization to higher dimensional manifolds. We conclude this section with the immediate inference that there exists a curvature-based Persistent Homology scheme for hypernetworks.

\subsection{The polyhedral surfaces case}

We begin by briefly recalling the basic ideas and notions of the classical Morse Theory. (For a full, yet succinct presentation we can not but recommend Milnor's classical \cite{Mil}). This will allow us to better explain the ideas residing behind Banchoff's polyhedral version.

Let $S^2$ be a smooth, closed surface  in $\mathbb{R}^3$ and let $\vec{\mathbf{v}}$ be a and arbitrary direction in $\mathbb{R}^3$ (i.e. a unit vector, or a point on the unit sphere $\mathbb{S}^2$). We define the {\it height function} $h$ as being the projection function of $\mathbb{R}^3$ on the line $l$ determined by $\vec{\mathbf{v}}$.  
A point $p \in S^2$ is called a {\it critical point} for $h$ if the tangent plane to $S^2$ at $p$ is perpendicular to $l$, otherwise it is called an {\it ordinary point}. 

To each critical point a numerical value is attached -- it's so called {\it index} (of $p$, with respect to the direction $l$), which is defined as follows: $i(p,l) = +1$ if $m$ is a local minimum or maximum, and $i(p,l) = -1$ if $m$ is a saddle point. This formal definition is motivated by the  following geometric observation: If $p$ is an ordinary point, then the tangent plane (to $S$) at $p$ is not ``horizontal'' (parallel to $l$), therefore it meets a ``small'' (infinitesimal) circle (on $S$) around $p$ in precisely two points. 
In contrast, the intersection of the tangent plane with such a circle at maximum or minimum point is void, whereas at a saddle point it will intersect an infinitesimal circle in four distinct points

Based on the observation above, Banchoff introduces  the following definition of the  index of the vertex $v$ of a polyhedral manifold (surface) $M^2$ in the following combinatorial manner:
\begin{equation}
\hspace*{-1cm}i(v,l) = 1 - \frac{1}{2}|\{T_N(v) \; \cap \; C_\varepsilon(v)\}|\,;
\end{equation}
where $T_N(v)$ is the plane through $v$ normal to $l$ and $C_\varepsilon(v)$ denotes an infinitesimal circle centered at $v$; that is
\begin{equation}
i(v,l) = 1 - \frac{1}{2}\big({\rm \#points\; in\; which\; the\; plane\; through\;} v\; 
\end{equation}
\[
{\rm \hspace*{2.5cm} perpendicular \; to}\;  l\; {\rm meets\; a\; ``small\; circle''\; about}\; v\; {\rm on}\; M^2\big).
\]

While quite intuitive, this definition is not precise enough for a notion for smooth surfaces, both because of the vagueness of the notion of ``small circle'', and because, in practice, it is rather difficult to determine the required  number of intersections on a general  surface, and for a general direction $l$ ($\vec{\mathbf{v}}$). On the other hand, its form is almost what one would request from a definition befitting polyhedral surfaces. 
To justify this assertion, let us first note that, for polyhedral surfaces,  the {\it star} $St(v)$ of a vertex $v$, i.e. the set of all simplices incident to $v$ (that is the edges and faces (including their edges and vertices containing $v$) plays the role of a ``small'' disk neighbourhood centered at $v$, while the {\it link} $Lk(v)$, i.e. the polygon representing the boundary of $St(v)$ represents the polyhedral analogue of a ``small circle'' around $v$. (For more details on this and other $PL$ Topology notions, see e.g. \cite{Hu}.)
Furthermore, observe also that a point $p$ is {\it ordinary} for the height function $h$ if the plane perpendicular to $l$ that passes through $p$  divides $St(p)$ into two pieces. Any interior point of  face of an edge is, therefore, an ordinary one for any direction {\it general} for the given polyhedral surface, i.e. such that $h(u) \neq h(v)$, for any two distinct vertices of $S^2$.  Moreover, given that $M^2$ has only a finite number of edges, it follows that the number of non-general directions is finite, thus our analysis is not limited by using general directions, given the fact that they are the rule, rather then one of the finite number of exceptional cases. Furthermore, this represents the precise polyhedral equivalent of a classical result, namely that almost any direction $\vec{\bf v} \in \mathbb{S}^2$, the associated height function has only a finite number of critical points, thus almost all directions (up to a set of zero measure)  are general (see \cite{Mil}).
Note also that, in contrast, vertices represent critical points of all of the possible types  for smooth surfaces. 
Moreover, while for smooth surfaces the only possible critical points are maxima, minima and non-degenerate saddle points (\cite{Mil}), 
on polyhedral surfaces 
{\it degenerate} critical points can also arise, like the so-called {\it monkey saddle} (see, e.g. \cite{Mil}). 

It is possible to improve this definition of the index to become purely combinatorial, by making the following observation:  One can count intersections with ``small disks'' even easier, because in the polyhedral context ``small disk'' has a precise meaning, namely the star of a vertex. 
It immediately follows from here that the number of times the plane through $v$ perpendicular to a triangle $T$ (with vertex $v$) meets ${\rm Lk}(v)$ is equal to $\# T$ in ${\rm St}(v)$, such that one of the vertices of $T$ lies above the plane and the other lies below. In such a case $v$ is called {\it the middle vertex} of $T$ for $l$.

The desired combinatorial definition of the index (at a vertex) is now easy to formulate:
\begin{equation}
i(v,l) = 1 - \frac{1}{2}\big({\rm \#{\it T} \; s.t. {\it p}\; is\; a\; middle\; for\; {\it l}}\big).
\end{equation}

The implications of  Banchoff's definition of the index a polyhedral manifold are both topological and geometrical, as they allow him to connect it both to the Euler characteristic and the (discrete) Gauss curvature (thus facilitating proofs both of Theorema Egregium and of the Gauss-Bonnet Theorem -- see \cite{Ban70}. The first main step towards proving these results is the fact that 
\begin{equation}  \label{eq:K(index)1}
K(v) = \frac{1}{2}\int_{\mathbb{S}^2}i(p,l)dA\,.
\end{equation}
(For details of the proof see \cite{Ban70}.)

Here the Gauss curvature at a point is the classical combinatorial one (going back seemingly to Descartes), namely
\begin{equation}
K(v) = 2\pi - \sum_{\alpha_i}\alpha_i\,;
\end{equation}
where the sum is taken over all angles $\alpha_i$ adjacent at $p$.

For combinatorial complexes, the triangles are equilateral, thus all angles equal $\pi$, thus the formula above becomes
\[
K(v) = \frac{\pi}{3}(6 - d_p)\,;
\]
where $d_p$ denotes the number of edges (and faces) incident at $p$.

This formula can be given a purely combinatorial flavor, by discarding the factor $\pi$ which is meaningless in a purely combinatorial context, such as that we'll adopt in the last section. Furthermore, can naturally obtained in a geometric context as well,  by the proper normalization of the area of the $\mathbb{S}^{n-1}$, an approach we adopt in the reminder of this section.
Thus  the formula for the curvature at a vertex becomes
\begin{equation} \label{eq:comb-K(v)}
K(v) = \frac{1}{3}(6 - d_p)\,.
\end{equation}


\subsection{The higher dimensional case}

The ideas presented above readily generalize to higher dimensions, and allow us to obtain simple connections between the sum of the indices of the vertices of a polyhedral manifold, its total curvature and the Euler characteristic of the manifold. More precisely, to define the index with respect to a general direction $l$, at a vertex $v$ of a polyhedral $n$-manifold $M^n$, one first has to introduce a characteristic function $I$, defined as follows:
\begin{equation}
I(C^k,v,l) = \left\{ 
\begin{array}{ll}
1\,;  &  \mbox{$p \in C^k$\, \; and $l(v) \geq l(u)$,  for\; all\; $u \in C^k$}\\
0\,;  & \mbox{else}\,. 
\end{array}
\right .
\end{equation}
where $C^k$ denotes a $k$-dimensional cell of $M^k$. Then the index of a vertex $v$ (relative to a general direction $l$) is defined as
\begin{equation}
i(v,l) = \sum_{k=0}^n(-1)^k\sum_{C^k \in M^n}I(C^k,v,l)\,.
\end{equation}

The first important result relates the index to the Euler characteristic. More precisely, we have the following polyhedral analogue of a classical result:

\begin{thm}[Banchoff \cite{Ban67}, Theorem 1] \label{thm:Banchoff-CPT} 
	Let $l$ be a general direction for the polyhedral manifold $M^n$ embedded in some $\mathbb{R}^N$. Then
	\begin{equation}
	\sum_{v\; {\rm vertex\; of\;} M^n}i(v,l) = \chi(M^n)\,.
	\end{equation}
\end{thm}

\begin{rem}
	Since a direction $l$ is general for $M^n$ if $l(u) \neq l(v)$, for any $u,v$ that are the end vertices of an edge of $M^n$, it easy to show that the set of general directions for a given polyhedral manifold is open and dense (\cite{Ban67}, Proposition 1).
\end{rem}

As already mentioned, there also is a close connection between the Euler characteristic and the curvature of the simplicial complex. The curvature in question is, as expected, the generalization of the classical {\it defect}, i.e. combinatorial one used in the 2-dimensional case. 
The curvature at a vertex $v$ is the vertex of a convex cell $C^r$ embedded in $\mathbb{R}^n$, it's curvature  is defined to be 
\begin{equation} \label{eq:Ban-CurvDef-HD}
K(v) = \frac{1}{\rm Area({\it \mathbb{S}^{r-1}})}\int_{\mathbb{S}^{r-1}}i(v,l)dA_{r-1}\,.
\end{equation}
and the total curvature of the complex as
\begin{equation}
K(M^n) = \sum_{v \in M^n}K(v) = \frac{1}{\rm Area({\it \mathbb{S}^{r-1}})}\int_{\mathbb{S}^{r-1}}\sum_{v \in M^n}i(v,l)dA_{r-1}\,.
\end{equation}

From this last formula and from Theorem \ref{thm:Banchoff-CPT} above one readily obtains the following generalization of the classical {\it Gauss-Bonnet Theorem}:

\begin{thm}[Banchoff \cite{Ban67}, Theorem 2] \label{thm-Ban-GB+}
	Let $M^n$ be a polyhedral manifold embedded in some $\mathbb{R}^N$. Then
	\begin{equation}
    K(M^n) = \chi(M^n)\,.
    \end{equation}
\end{thm}

Having thus obtained a proper, geometrically intuitive Morse theory, thus automatically a filtration method for Persistent Homology, we are still faced with the question regarding the practical feasibility of our approach. In point of fact, while the index and curvature of a vertex are easily computable, the definitions above depend on a specific embedding of the complex in some higher dimensional Euclidean space. 
%
%
In truth, while for simplicial complexes this is easier than for general polyhedral ones (\cite{Ban67}, \cite{Ban70}, \cite{Ban83}), this approach, without further refinement is not truly feasible for complexes of dimensions higher than 3, since the intuitive aspect is lost, at least as far as practical computations are concerned. However, Banchoff has also shown that, akin to the concept for smooth manifolds, the curvature of a vertex is {\it intrinsic}, i.e. it does not depend on the specific embedding. To be able to technically formulate this result and to understand how to compute the curvature of a vertex in an intrinsic manner, 
we first need to recall the notion of {\it normalized exterior angle}: Given a convex cell $C^k \subset \mathbb{R}^n$, and a vertex $v$ of $C^k$, the normal exterior angle of $C^k$ at $v$ is defined as:
\begin{equation}
E(C^k,v) = \frac{\rm Area(set\; of\; normals\; to\; support\;  hyperplanes\; at\; {\it v})}{\rm Area({\it \mathbb{S}^{k-1}})}\,\;
\end{equation}
that is
\begin{equation}
E(C^k,v) = \frac{1}{\rm Area({\it \mathbb{S}^{k-1}})}\int_{\mathbb{S}^{k-1}}I(C^r,v,l)dA_{k-1}\,.
\end{equation}
This definition is independent (\cite{Ban67}, Lemma 2) on the dimension of poyhedral complex of which $C^k$ is a cell: If $C^k \subset C^k \subset C^n$, then 
\begin{equation}
E(C^k,v) = \frac{1}{\rm Area({\it \mathbb{S}^{n-1}})}\int_{\mathbb{S}^{n-1}}I(C^k,v,l)dA_{n-1}\,.
\end{equation}

The {\it intrinsic definition of curvature} at a vertex is then given as follows:
\begin{equation}
\kappa(v) = \sum_{i = 0}^n\sum_{C^k \in M^n}E(C^k,v)\,.
\end{equation}

With this definition one can formulate and prove the polyhedral analogue of the classical {\it Theorema Egregium}:

\begin{thm}[Banchoff \cite{Ban67}, Theorem 3] \label{thm: Ban-CurvThm-HighD}
	$K(v) = \kappa(v)$, therefore $K(v)$ is intrinsic. 
\end{thm}

Thus Banchoff's index-based definition of curvature coincides to the classical defect one, certainly for 2-dimensional {\it (pseudo-)manifolds} (i.e. combinatorial manifolds with possible singularities) and manifolds with boundary, but also for more general cases (see \cite{Ban67} Remarks 5 and 6, respectively). Thus for the unweighted hypernetworks, the purely combinatorial form of curvature (\ref{eq:comb-K(v)}) is the natural one to consider, especially in applications. Moreover, given that one looks for a {\em topological} filtration, a this version can be employed as well in the general case. We shall dwell upon the combinatorial curvature in more detail in the sequel. 
In fact, there is not necessary for the complex $M^n$ to be embedded in $\mathbb{R}^n$, neither for the proof of the essential Theorem \ref{thm:Banchoff-CPT}, nor for the derivation of the curvature related theorems. Moreover, the restriction to simplicial complexes is not necessary 
as the results can be extended to {\em abstract} cell complexes and it is possible to replace the embeddings by much less restrictive {\em local isometric mappings} into $\mathbb{R}^n$ (\cite{Ban67}, Remark 7) -- See also the discussion at the end of the next section.) 

More precisely, one considers finite, connected complexes composed of bounded convex cell, such that the face inclusions for each cell are given by isometric maps. The directions (thence height functions) are cell-wise linear functions and will be called generic if they are homeomorphisms when restricted to each 1-dimensional cell. These modifications suffice to prove the critical point theorem Theorem \ref{thm:Banchoff-CPT}  in this generalized setting.
To attain the analogue of Theorem \ref{thm: Ban-CurvThm-HighD} one has to consider local isometric mappings $f$ of $K^n$ (equipped with its intrinsic metric) into $\mathbb{R}^N$. Then the  definition of curvature at a vertex becomes
\begin{equation}
K(v) = \frac{1}{\rm Area({\it \mathbb{S}^{r-1}})}\int_{\mathbb{S}^{r-1}}i(v,l \circ f)dA_{r-1}\,.
\end{equation}
Thus, the main obstruction in the practical application of the Banchoff Morse Theory to the Topological Data Analysis of hypernetworks, namely the necessity of using the restrictive isometric embeddings, is removed, and one can employ it by computing the index of each vertex, via the simple combinatorial curvature mentioned above.

Moreover, it turns out that even the condition that $l$ be general w.r.t. $M^n$ is superfluous, once a more general definition of the index is achieved  see \cite{Ban67}, Remark 7 ff..
We begin by introducing a generalized indicator function, as follows:
\begin{equation}
I(C^k,C^j,l) = 	
\left\{ \begin{array}{ll}
1; & C^j \subset \partial C^k, \geq, \forall x \in C^j, y \in C^k\,;\\
0\,; & {\rm else}\,.
\end{array}
\right.
\end{equation}
(For vertices $v = C^0$ this definition reduces to the one already considered above.)

Using this relaxed notion of index, one can obtain the following result:

\begin{prop}[Banchoff \cite{Ban67}, Corollary 4]
	Let $M^n$ be a convex cell complex. Then %
	\[\chi(M^n) = \sum_k^n(-1)^k\sum_{C^k \in M^n}\chi(St(C^k))\,;\]
	where $St(C^k)$ denotes the {\it open star} of $C^k$, $St(C^k) = \{C^m \in M^n | C^k < C^m\}$\,.
\end{prop}

Thus it is possible to obtain global topological information from the local one in every dimension, that can be easily read from the given network, in its original form.

We can summarize the discussion above, that prescribes to hypernetworks (and multiplex networks), viewed as polyhedral complexes in $\mathbf{R}^N$, a natural Morse function, as the following practical result:

\begin{thm}
	There exists a curvature-based Persistent Homology scheme for hypernetworks. 
\end{thm}


\section{The Connection between Forman and Banchoff Morse Theories}

As we already noted above, as a {\em bona fide} Morse Theory, Banchoff's polyhedral version allows for a filtering scheme for hypernetworks. Its intuitiveness makes it simple and attractive, especially in the low dimensional case. The question still remains whether it is possible to formulate this as a precise algorithm whose properties can be rigorously analyzed? There exist at least one way to arrive to a positive answer. However, before presenting it, we must first present the connection between Banchoff's Morse Theory and yet another Discrete Morse Theory, due to Forman \cite{Fo98}.

Forman's Morse theory applies to a class of geometric objects that is more general than polyhedral manifolds (and that contains it), namely that of ({\it regular}) {\it CW complexes}.
As such, this type of Morse Theory is both more ``discrete'' and, in consequence, less intuitive than Banchoff's version. 
Since this approach is not the one we adopt here and, moreover, it is neither specific to, nor especially descriptive for polyehdral manifolds, which represent our model for hypernetworks, we do not detail it here, but we rather direct the reader to Forman's original paper \cite{Fo98}. (See also \cite{KSRS})

While these two discrete notions on Morse Theory are, as mentioned above, quite dissimilar, they are, surprisingly, interconnected, as discovered by E. Bloch \cite{Bloch13}. More precisely, we have the following result 

\begin{thm} [Bloch \cite{Bloch13}, 2013] \label{thm:Bloch-Fvs-B}
	Let $X$ be a finite regular $CW$ complex, and let $f$ be a discrete Forman-Morse function on $X$. For any sufficiently large $m \in \mathbb{N}$, and for any unit
	vector $\vec{v} \in \mathbb{S}^{m-1}$, there is a polyhedral embedding of the barycentric subdivision of $X$ in $\mathbb{R}^m$ such that a cell in $X$ is Forman-critical with respect to $f$ if and only if its barycenter is Banchoff-critical with respect to projection onto
	the line spanned by $\vec{v}$.
\end{thm}

\begin{rem}
	The necessity of taking barycentric subdivisions follows first of all from the fact that while Banchoff's Morse Theory determines the critical {\it vertices}, Forman's version asserts the criticality of {\it cells}. Thus it is imperative to be able to replace, when using Banchoff's approach, each cell with a corresponding vertex. 
	
	This subdivision process also proves to be actually an advantage, since embedding of simplicial complexes is standard \cite{Fl}, \cite{vK} and far easier than the one of general $CW$ complexes.
	
	Furthermore, while the function assigning to each face of a simplex triangle its dimension is a Forman-Morse function, and, moreover, each face is critical with respect to this discrete Morse function,  the projection onto a general line (direction) in $\mathbb{R}^N$ space assigns to any point in the interior of an edge a lower hight than precisely one of its vertices, and greater than the other its vertices \cite{Fo98a}, thus the barycenter of an edge is Banchoff-Morse ordinary. Therefore, to ensure the the Forman-Morse criticality of simplices to coincide with the Banchoff-Morse one, it is necessary to perform the barycentric subdivision before the embedding into Euclidean space.
	
	Finally, let us note in this context that the minimal embedding dimension is not stable under barycentric subdivisions (in fact, it can alternatively increase and decrease under successive subdivisions) -- see \cite{Cairns}, Theorem IV. Even though the theorem assures 
	 the existence  of a large enough $N$, in practice it is important to find the minimal such $N$, especially if visualization is also required. This further strengthens the need for taking the barycentric subdivision prior to the embedding.
\end{rem}

We can therefore pass from a hypernetwork realized in some Euclidean space and the discrete Morse function associated to its vertices, take advantage of techniques, algorithms  and results in \cite{KSRS}, and thus obtain the following 

\begin{thm} \label{thm:BtF}
	There exists a Discrete Morse Theory-based Persistent Homology scheme for hypernetworks that achieves a close to theoretical minimum number of critical simplices. 
\end{thm}

It remains to ascertain which of the two versions of discrete Morse Theory is more advantageous in practice. On the side of the Banchoff approach is its intuitiveness, at least in dimension $\leq 3$. While the index and curvature computation retain their simplicity, the geometric intuition and the visualization readiness are lost in dimension higher than 3. This represents a relative weakness of this method. However, we should also take into account that, in practice, no true experiments (at least not on medium to large scale data) were performed with any Persistent Homology method in dimension higher than 3 (this holds, in particular, for the experiments in \cite{KSRS}, \cite{RVRS}). Still, the Forman-Morse approach appears to be, while less intuitive, also more readily applicable to raw, abstract data. Here we should however note that when passing to the Banchoff approach via barycentric subdivision, the number of simplices of the resulting complex increases drastically as compared to the one of the original one, with clear computational implications. Thus applying the Forman-Morse function to a hypernetwork embedded as a simplicial/polyhedral complex in Euclidean space will be computationally more costly than ideal, due to the fact that one has to consider all the vertices, and not just those corresponding to the ``hidden'' cells of the cell complex optimal  vis-a-vis the Forman-Morse approach. Indeed, even if such vertices will have precisely $n+1$ neighbors ($n$ being the dimensionality of the complex), they are by no means the only ones satisfying this property.

Before concluding this section, let us note that there is yet another advantage in adopting the Banchoff-Morse theory, rather than Forman's one (at least in feasible dimensions): Forman-Morse based filtration devised in \cite{KSRS} that resides behind Theorem \ref{thm:BtF} is applicable solely to {\em unweighted} hypernetworks. However, real-life networks and hypernetworks are vertex- and edge-weighted networks. By metricizing the hypernetwork/complex (e.g. by using, as suggested in \cite{WJS1}, the {\it path degree metric} \cite{DK}), one obtains a polyhedral realization in some $\mathbb{R}^N$ of the hypernetwork, a realization to which the Banchoff approach can be applied. Moreover, this can be done in a manner that takes into account the weights, by considering the extrinsic curvatures determined by the weight-induces lengths of the edges.


\section{Bloch's Stratified Index}

The main obstruction to the full applicability of the Banchoff's Morse Theory to the Persistent Homology of hypernetworks resides in the fact that it holds only for $PL$ manifolds and the somewhat more general convex cell complexes, thus the essential equivalence between Banchoff's and Forman's Morse holds only for this type of spaces (a fact that is somehow omitted in Bloch's original formulation of Theorem \ref{thm:Bloch-Fvs-B}). However, Bloch also developed a {\it stratified Morse Theory} that extends Banchoff's ideas and work to general simplicial complexes. 
Both for the sake  af completeness and because these ideas and results are far less familiar to most readers than Banchoff's work, we present them succinctly here. However, for most of the details, including examples, further motivation and proofs we refer the reader to Bloch's original papers \cite{Bloch04} and \cite{Bloch98}.

In the following, we shall always denote by $K = K^n$ a finite $n$-dimensional simplicial complex embedded in $\mathbb{R}^N$, for some large enough $N$; by $\sigma = \sigma^n$ an $n$-dimensional simplex of $K$, and by $\tau = \tau^k$ a $k$-dimensional face of $\sigma$, that is $\tau < \sigma \in K$.

To define the {\it stratified Morse Theory} that generalizes Banchoff approach as developed by Bloch \cite{Bloch04}, we must first introduce a number of notations and definitions.

\begin{defn}
We denote by $T_i$ the open cone on $i$ points, $T_0 = {0}$, and define 
\begin{equation}
P_{n,i} = T_i \times \mathbb{R}^{n-1}.
\end{equation}
Furthermore, the set $\{*\} \times P_{n,i}$, where $\{*\}$ denotes the cone point of $T_i$, is called the {\it apex set} of $T_i$.
\end{defn}

\begin{rem}
If $|K|$ is a manifold, any point has a neighborhood homeomorphic to $P_{n,2} \simeq (-1,+1) \times \mathbb{R}^{n-1} \simeq \mathbb{R}^n$, whereas if it is a manifold with boundary, it has a neighborhood as above if it is an interior point, or homeomorphic to $P_{n,1} \simeq [0,+1) \times \mathbb{R}^{n-1}$. Also, while not every point in $|K|$ has a neighborhood homeomorphic to $P_{n,2}$ or $P_{n,1}$, each point in $|K| \setminus |K^{(n-2)}|$ has a neighborhood homeomorphic to some $P_{n,r}$, where $r$ varies depending an the point.
\end{rem}

\begin{defn}
	Let $0 \leq r$ be a natural number. We define the $C_r^n(K) \subset K$ as follows:
	\begin{equation}
C_r^n(K) =	\left\{ \begin{array}{ll}
		  \{x \in |K|\,|\, \exists V \in \mathcal{N}(x), {\rm s.t.} V \stackrel{h}{\simeq} P_{n,r}\,, h(x) = * \}; & r \neq 2\,;\\
		  |K| \setminus \bigcup_{r \neq 2}C_r^n\,; & r = 2\,.
		  \end{array}
		  \right.
	\end{equation}
\end{defn}

\begin{rem}
	\begin{enumerate}
		\item $\{C_r^n\}_r$ represents a cover of $K$ with disjoint sets. Moreover, for $r \neq 2$, $C_r^n$ is an $(n-1)$-manifold, and it is the union of open simplices of $K$. 
		
		\item If $|K|$ is an $n$-manifold, then $C_2^n(K) = |K|$ and  $C_f^n(K) = \emptyset$ if $r \neq 2$; and 
		        if $|K|$ is an $n$-manifold with boundary, then $C_2^n(K) = |K| \setminus \partial|K|$  and  $C_f^n(K) = \emptyset$ if $r \neq 2$.
		
		\item $|K| \simeq |L| \rightarrow C_r^n(K) \simeq C_r^n(L)$.
	\end{enumerate}
\end{rem}

We next introduce a type of ranking for the simplices of $K$:

\begin{defn} 
Let $K$ be as above and $\sigma$ a simplex of $K$. We define:  $T_n(\sigma) = r/2$, where $\sigma \in C_r^n(K)$ for a unique integer $r$.
\end{defn}

\begin{rem}
	\begin{enumerate}
		\item If $\sigma = \sigma^n \in K^n$, then $T_n(\sigma) = 1$, and if $\tau = \tau^{n-1} \in K^n$, then $T_n(\tau) = \frac{1}{2}|\{{\rm vertices\; in\;} Lk(\tau)\}|$.
		
		\item $T_n(\sigma)$ is not locally computable, since it is depends on the dimension of $K$, which is not a local quantity. However, if $K^n$ is {\it purely $n$-dimensional},  (in particular when $K^n$ is an $n$-dimensional $PL$ manifold), then $T_n(\sigma)$ is locally computable. (Recall that $K^n$ is called  {\it purely $n$-dimensional} if for for any simplex $\tau^k \in K^n$, there exists $\sigma^n \in K^n$ such that $\tau^k < \sigma^n$.)
	\end{enumerate}
\end{rem}

Bloch's approach rests on the notion of {\it generalized angle defect} as defined below:

\begin{defn}
	The {\it generalized angle defect} at $\eta^i$, is denoted as $D_n(\eta^i)$, and is defined as follows:
	
\begin{equation}
D_n(\eta^i) = T_n(\eta^i)	- \sum_{\eta^i \in \sigma^n}\alpha(\eta^i, \sigma^n)\,;
\end{equation}
where $\alpha(\eta^i, \sigma^n)$ denotes the {\it solid angle in $\sigma^n$ along $\eta^i$}, normalized such that $\alpha(\eta^i, \sigma^n) \in [0,1]$ (i.e. such that ${\rm Vol}_{n-1}(\mathbb{S}^{n-1}) = 1$, for all $n$).	
\end{defn}

Note that, in contrast to the standard approach to discrete curvature of polyhedral manifolds, where curvature is concentrated solely at vertices, Bloch's definition assigns curvature to all the simplices (even though it is non-zero only at simplices of codimension $\geq 2$ and the angle sum changes only along codimension 1 singularities). -- See also \cite{Ban83} for a similar approach.

We can now introduce the Bloch's definition of {\it stratified Euler characteristic}:

\begin{defn}
The {\it stratified Euler characteristic} of $K$ is defined as
\begin{equation}
	\chi^s(K) = \sum_{\eta \in K}T_n(\eta)(-1)^{{\rm dim} \eta}
\end{equation}	
\end{defn}

\begin{rem}
	\begin{enumerate}
		\item $\chi^s(K)$ is not, in general, an integer. However, $\chi^s(K) = p/2\,,\: p \in \mathbb{Z}$\,.
		
		\item $\chi^s(K)$ is not a homotopy type invariant (in contrast with the classical Euler characteristic), however it depends solely on $K$' up to homeomorphism. In particular, it does not depend on the specific triangulation of $|K|$ that renders $K$.
	\end{enumerate}
\end{rem}

Bloch's generalization of Banchoff's notion of a Morse Theory for $PL$ (polyhedral) manifolds rests (naturally) on a fitting concept of a {\it stratified Morse index}. However, before we can present it, we have first to bring a number of preparatory notions and definitions.

We begin by introducing some notations:
\begin{itemize}
	\item For any simplex $\eta^i \subset \mathbb{R}^N$, we denote by $V(\eta^i)$ the $i$-dimensional vector space of $\mathbb{R}^N$ parallel to the $i$-plane spanned by $\eta$.
	
	\item If $U$ is a vector subspace of $\mathbb{R}^N$, we denote by $h_U: \mathbb{R}^N \rightarrow U$, the orthogonal projection on $U$. We put $h_{\vec{\mathbf v}} = h_{V(\vec{\mathbf v})}$, for any vector $\vec{\mathbf v}  \in \mathbb{R}^N$. Here we view $V(\vec{\mathbf v})$ as a copy of $\mathbb{R}$, thus we also regard $h_{\vec{\mathbf v}}(x)$, $x \in \mathbb{R}^N$, as a real number, rather than a vector. Note also that, for any unit vector $\vec{\mathbf v}$, $h_{\vec{\mathbf v}}(x) = x \cdot \vec{\mathbf v}$. 
\end{itemize}

In order to define Morse functions on $K$ we have to use projections of the form $h_{\vec{\mathbf v}}:\mathbb{R}^N \rightarrow V(\vec{\mathbf v})$. However, to do so, one has first to discard the ``bad'' vectors, which are formally defined as follows:

\begin{defn}
	Let $K$ be as above and let $\vec{\mathbf v} \in \mathbb{S}^{m-1}$. Then $\vec{\mathbf v}$ is called and {\it allowable vector} with respect to $K$ if the following conditions hold:
	\begin{enumerate}
		\item $h_U(\vec{\mathbf v}) \neq \vec{0}$\,;
		
		\item $h_U(\vec{\mathbf v}) \nsubset V(\eta^{n-1})$, for any $\eta^{n-1} < \sigma^n$;
	\end{enumerate}
where $U$ is chosen, without loss of generality, as $U = V(\sigma)$, for any $n$-simplex $\sigma \in K$.
\end{defn}

Fortunately, one can discard the ``bad'' (non-allowable) vectors due to the fact (cf. \cite{Bloch04}, Lemma 3.2) that 

\begin{enumerate}
	\item The set of allowable vectors with respect to $K$ is an open, dense subset of $\mathbb{S}^{m-1}$;
	
	\item  The set of non-allowable vectors has zero measure.
\end{enumerate}

\begin{defn}
Let $K, \sigma = \sigma^n$ and $\vec{\mathbf v}$ allowable w.r.t. $K$ be as above, and suppose, without loss of generality, that $U = V(\sigma)$.
the unit vector
\begin{equation}
\vec{\mathbf v}_{\sigma} = \frac{h_U(\vec{\mathbf v})}{||h_U(\vec{\mathbf v})||}
\end{equation}
is called {\it allowable w.r.t.} $\sigma$.
\end{defn}

Before being able to define the stratified index, we need to introduce a number of other preparatory notions.

\begin{defn}
	Let $K = K^n \subset \mathbb{R}^N$ be a simplicial complex;  let $\sigma^n \in K, \sigma^n = <a_0,a_1,\ldots,\ldots,a_{n-1};a_n>$; $\tau^{n-1}  < \sigma^n, \tau^{n-1} = <a_0,a_1,\ldots,a_{n-1}>$; and let $\vec{\mathbf v}$ be an allowable vector w.r.t. $K$.
	
	Denote: $x_i = a_i - a_0, i = 1, \ldots, n - 1$; $y = a_n - a_0$. 
	
	We define 	
	\begin{equation}
	t(\tau^{n-1},\sigma^n,\mathbb{R}^N,\vec{\mathbf v}) =	
	\left\{ \begin{array}{ll}
	1; & D > 0\,;\\
	0\,; & {\rm else}\,.
	\end{array}
	\right.
	\end{equation}
	where 
	
	\begin{equation}
	D = \left| \begin{array}{cccc}
		x_1 \cdot x_1 & \ldots & x_{n-1} \cdot x_1 & h_{\vec{\mathbf v}}(x_1) \\
		\vdots &  & \vdots & \vdots  \\
		x_1 \cdot x_{n-1}  & \ldots & x_n \cdot x_{n-1} & h_{\vec{\mathbf v}}(x_{n-1})  \\
		x_1 \cdot y  & \ldots & x_1 \cdot y & h_{\vec{\mathbf v}}(y)
	\end{array}
	\right| \,.
		\end{equation}
	
\end{defn}

\begin{rem}
	\[t(\tau^{n-1},\sigma^n,\mathbb{R}^N,\vec{\mathbf v}) = 1 \Longleftrightarrow \vec{\mathbf v}_{\sigma^n} {\rm points\; accross}\; \tau^{n-1}\; {\rm in\; the\; direction\; of}\; \sigma^n;\]
	\[\Longleftrightarrow  h_{\vec{\mathbf v}}(y) > h_{\vec{\mathbf v}}(x),  \forall x \in {\rm int}\tau^{n-1}   {\rm and }\; \forall y \in {\rm int}\sigma^n,\; {\rm such\; that}\; \overrightarrow{y-x} \perp V(\tau^{n-1})\,.\]
\end{rem}

\begin{defn}
Let $K, \sigma^n \in K$ as above, let $\eta^j < \sigma^n, 0 \leq j \leq n - 2$, and let $\vec{\mathbf v}$ be an allowable vector w.r.t. $K$. We define
\begin{equation}
g(\eta^j, \sigma^n, \mathbb{R}^N,\vec{\mathbf v}) = \prod_{k=1}^{n-j}t(\tau^{k},\sigma^n,\mathbb{R}^N,\vec{\mathbf v}) + \prod_{k=1}^{n-j}t(\tau^{k},\sigma^n,\mathbb{R}^N,-\vec{\mathbf v})\,;
\end{equation}
where $\eta^j = \bigcap_{k=1}^{n-j}\tau_k$\,.
\end{defn}

We can no introduce the definition of the {\it stratified index}:

\begin{defn}
et $K, \sigma^n \in K$ as above, let $\eta^j < \sigma^n, 0 \leq j \leq n - 2$, and let $\vec{\mathbf v}$ be an allowable vector w.r.t. $K$. The {\it stratified index of $\eta^j$ relative to $\vec{\mathbf v}$} is defined as follows
\begin{equation}
i^s(\eta^j,\vec{\mathbf v}) = T_n(\eta^j) - \frac{1}{2}\sum_{\sigma^n > \eta^j}g(\eta^j, \sigma^n, \mathbb{R}^N,\vec{\mathbf v})\,;
\end{equation}	
where the sum is taken over all $\sigma^n \in K$, such that $\eta^j < \sigma^n$.	
\end{defn}

We can now bring the first of Bloch's results essential in a stratified Morse Theory:

\begin{thm}[Stratified Critical Point Theorem; \cite{Bloch04}]
Let $K = K^n \subset \mathbb{R}^N$ be a simplicial complex and let $\vec{\mathbf v}$ be an allowable vector w.r.t. $K$. 
Then 
\begin{equation}
\sum_{\eta^j \in K,\: 0 \leq j \leq n-1}(-1)^ji^s(\eta^j,\vec{\mathbf v}) = \chi^s(K)\,.
\end{equation}

\end{thm}

The second essential ingredient in the applying Bloch's generalized (stratified) index to Morse Theory for general simplicial complexes is the stratified version of Theorema Egregium below:

Given that the set of allowable vectors is a dense, open set of $\mathbb{S}^{m-1}$, and its complement is a set of zero measure, we may consider the integral of 

\begin{thm}[Stratified Theorema Egregium; \cite{Bloch04}]
Let $K = K^n \subset \mathbb{R}^N$ be a simplicial complex and let $\eta^j \in K, 0 \leq j \leq n - 2$. 
Then
\begin{equation}
\int_{\mathbb{S}^{m-1}}i^s(\eta^j,\vec{\mathbf v})d{\rm Area}^{m-1} = D_n(\eta^j)\,.
\end{equation}
%
\end{thm}
(Note that the integral above is well defined due to the fact that the set of allowable vectors is a dense, open set of $\mathbb{S}^{m-1}$, and its complement is a set of zero measure.) 

These theorems render the desired Morse Theory for general simplicial complexes, that can then be applied for the Persistent Homology of hypernetworks, as well as multiplex networks, not only to the much more ``tame'' {\it clique complex} as in \cite{KSRS}). To this end, note that even if it is less intuitive than Banchoff's version, 
the index computation is essentially combinatorial in nature and as such applicable to any dimensions. 
Let us also note here that using the stratified Morse Theory, one can compute the Persistent Homology of quite general structures. Indeed, even though we show in \cite{SS} that any hypernetwork can be construed, in a canonical manner, as a simplicial complex, the transformation applied to this end to the network, while simplifying the picture, also obliterates much relevant geometric information. Moreover, while Bloch's Morse Theory requires -- as Banchoff's does -- for the complex to be embedded in some Euclidean space, this condition can also be significantly relaxed, given the fact that the notion of solid angle can be consistently defined for more general {\it simplex-wise embedded} complexes. (Recall that a simplex-wise embedding of a simplicial complex is a pair $(K,f)$, where $K$ is an $n$-dimensional simplicial complex and $f:K^(0) \rightarrow \mathbb{R}^N$ (for some $N$), such that if $a_0,\ldots,a_n$ are the vertices of a simplex $\sigma^n \in K$, then 
$f(a)_0,\ldots,f(a)_n$ are affinely independent in $R^N$. Note that $f$ is not necessarily and immersion, hence, a fortiori, not an embedding.) This extension to non-embedded complexes has significant practical implications, since it effectively allows for the computation of Persistent Homology in a manner that is as efficient yet more simple that the one based on Forman's Morse Theory. However, we have to also note again that Bloch's approach holds only for simplicial complexes, ever though quite general ones, and does not extend to general cell complexes. 

We shall, however, further simplify in the following section the computation of defects and express them in terms of Forman's Ricci curvature \cite{Fo}, thus both extending to hypernetworks as well as simplifying the scheme developed in \cite{RVRS}.


\section{ The Connection Between Defect and Forman-Ricci Curvature}

It is only natural to ask the whether there is a connection between the Persistent Homology approaches discussed above and Forman's Ricci curvature, and thus to the other topological aspects and applications of this type of curvature. This question is even more pertinent in light of the connection between the Banchoff-Morse index and curvature on the one hand, and the recent strong numerical correlations between the  Persistent Homology results given by the Forman discrete Persistent Homology and the based (alter alia) using Forman's Ricci curvature, that were recently observed on model and real life networks \cite{RVRS}. 

We shall show that the answer is, perhaps as expected, positive. 
To this end we make appeal to the following simple formula below  (see \cite{Fo}): 
\begin{equation} \label{eq:Forman-e-Comb}
{\rm Ric}(e) = \sharp\{2-{\rm cells}\; f > e\} + 2 - \sharp\{{\rm parallel\; neighbors\; of}\; e \}.
\end{equation}
which is the natural one to consider given the combinatorial type of curvature to which we want to relate. 
Recall  that edges $e$ and $\hat{e}$ are called {\it parallel}, if either $e$ and $\hat{e}$ belong the boundary of a common 2-face (2-cell), or have a common vertex, but not both these  incidences occur simultaneously.

In the case of {\it purely 2-dimensional} simplicial complexes (i.e. for which all triangles have faces ``plumbed in''), where none of the sides of the triangular faces $t$ adjacent to $e$ are parallel to it, the formula above reduces to
\begin{equation} \label{eq:Forman-e-Comb}
{\rm Ric}_F(e) = |\{t  > e\}| + 2 - (d_u + d_v - 2 - 2|\{t  > e\}|) = 4 - 3|t| - (d_u + d_v) \,;
\end{equation}
where $u$ and $v$ are the vertices of $e$. (In the case of general simplicial complexes, one can only ascertain that $4 - 3|t| - (d_u + d_v)$ represents the upper bound for ${\rm Ric}(e)$.) For regular complexes, i.e. such that $d_v \equiv d$, for some natural number $d$, the formula above becomes
\begin{equation} \label{eq:Forman-e-Comb1}
{\rm Ric}_F(e) 
= 4 - 3|\{t  > e\}| - 2d \,;
\end{equation}
i.e. reduces to a counting formula for the triangles containing the edge $e$.

While, as mentioned in the introduction, we do not expand here on the notion of Ricci curvature, but rather refer the reader to our previous papers on its applications to Complex Networks cited above and, of course, to Forman's original work \cite{Fo}, we wish to emphasize that, since it is a function on edges, it is determined solely by the 1- and 2-faces of the network (complex).

Furthermore, being a vertex measure, combinatorial curvature is, in fact, a scalar curvature.\footnote{For a brief overview of curvatures in Riemannian Geometry and their discretizations, see \cite{Sa-Metric1}.} It is therefore proper to compute it as such, using the Forman-Ricci curvature of the edges incident to a vertex, namely
\begin{equation}
{\rm scal}_F(v) = \frac{1}{d_v}(F(e_1) + \cdots F(e_d))\,.
\end{equation}
%
Since in a $PL$ (pseudo-)manifold each edge is common to precisely two 2-faces, and since these faces are all triangles, there are no parallel edges except the edges incident to the end points of $v$, it follows from Formula (\ref{eq:Forman-e-Comb}) that the expression of ${\rm scal}_F(v)$ becomes
\[{\rm scal}_F(v) = \frac{1}{d}[(4 - d -1 - d_1 - 1)  + \cdots (4 - d - 1 - d_d) - 1] = \frac{1}{d}\sum_1^d(2 - d - d_d)\]
\[\hspace*{-5cm} = 2 - d + \frac{1}{d}(d_1 + \cdots + d_d)\,,\]    
where $d = d(v)$ denotes the degree of the vertex $v$.
Thus the Forman-scalar curvature at a vertex equals %
\[2 - d - d_{\rm mean}\,,\]
where $d_m$ denotes the mean vertex degree of the complex $X$.
For (vertex) regular complexes, this equals $2 - 2d$. Furthermore, for ``almost regular'' hypernetworks/complexes, that is networks where almost all vertex degrees are equal  to the mean degree, such as those used in Graphics (after mesh improvement), where the degree is determined by the background topology of the network, we obtain

\[{\rm scal}_F(v) \approx 2 - 2d\,,\]
or
\[{\rm scal}_F(v) \approx 2 - 2d_{\rm mean}\,,\]
or yet again
\begin{equation} \label{eq:Fmean}
{\rm scal}_{F,{\rm mean}} = 2 - 2d_{\rm mean} = 2(1 - d_{\rm mean})\,.
\end{equation}


Unfortunately, a simple formula akin to (\ref{eq:Fmean}) connecting the Forman scalar curvature of a vertex solely to its degree is not possible, due to the fact that the former is defined in terms of the Ricci curvatures of the edges incident to the given vertex, thus as a function of the degrees of the vertices adjacent to $v$. The best one can obtain is, for example
\begin{equation} \label{eq:Forman-vs-Combi1}
{\rm scal}_F(v) = 3K(v) - 4 -  d_{\rm mean}\,. 
\end{equation}
While less elegant than perhaps wished for, this formula demonstrates that one can express the combinatorial curvature, thus the index of a vertex in an embedded $PL$ complex/hypernetwork in terms of its Forman scalar curvature and therefore, produce a Persistent Homology scheme in terms of Forman curvature. 
For 2-dimensional (closed) $PL$-surfaces, the mean curvature at each vertex is prescribed by the surface's topology. Indeed 
\[K_{\rm mean} = \frac{1}{|V|}\sum_{v \in V}K(v)\,. \]
On the other hand, we have by the discrete Gauss-Bonnet Theorem \ref{thm-Ban-GB+} above that 
\[ \sum_{v \in V}K(v) = \chi(X) = 2 - 2g\,,\]
where $g$ is the genus of $X$, we have
\[
K_{\rm mean} = \frac{1}{|V|}(2 - 2g)\,.
\]
Since
\[
K_{\rm mean} = 2 - \frac{d_{\rm mean}}{3}\,.
\]
From the last two equalities we obtain that
\[
d_{\rm mean} = 3\left[ \frac{1}{|V|}(2 - 2g) - 2\,\right]\,.
\]
It follows from here and from Formula (\ref{eq:Forman-vs-Combi1}) that, for $PL$-surfaces the Forman scalar curvature at a vertex is connected to the combinatorial one via the following nice formula that displays the important role of the background topology in the relationship between defect and Forman's scalar curvature: 
\begin{equation} \label{eq:Forman-vs-Combi2}
{\rm scal}_F(v) = 3\left[K(v) -  \frac{1}{|V|}(2 - 2g) + 2\,\right]\,.
\end{equation}
%

\begin{rem}
	Note that while a 0-dimensional analogue of  (\ref{eq:Forman-e-Comb}) is also available and it is tempting to implement it directly, without making appeal to the edge (Ricci) curvature, it is a useless essay, since in this degenerate case it produces, for PL surfaces, constant 0 curvature for each vertex, $F(v) = 0$, for all vertices $v \in X$. (Recall that vertices are not ``parents'' for any lower dimensional faces, thus the second term of  (\ref{eq:Forman-e-Comb})'s analogue is zero.)  
\end{rem}

Given that ${\rm scal}_F(v)$ is the mean of the Ricci-Forman curvatures of a edges adjacent to $v$, and these curvatures, as shown by conform Formula (\ref{eq:Forman-e-Comb})  --  depend in turn not just on the degree of $v$, but also on the degrees of the vertices adjacent to $v$, one can not hope, in the general case, to obtain tighter estimates that do not depend on the combinatorics of the 1-start of $v$ or, as above, on the global topology of the complex.


\section{Discussion and  Future Work}

We have summarized above the research regarding geometric Morse Theory for simplicial and more general cell complexes and showed that this represents a viable alternative, much more intuitive and easily implementable, at least in low dimensions, that Forman's better known -- and by now widely employed -- combinatorial approach. Given our interpretation of hypernetworks and multiplex networks \cite{SW19}, \cite{SS} and, indeed, the role of polyhedral complexes and their curvatures in the understanding of networks in general \cite{WJS}, \cite{SSGLSJ}, this amounts to existence of geometric schemes for Persistent Homology computation of networks and their generalization. Furthermore, we have shown the correlation between Forman's combinatorial Ricci curvature and the combinatorial defect at vertices, therefore between the former one and the index, thus showing that Forman's curvature can be employed for the Persistent Homology of (hyper-)networks. Given the relationship between Forman's and Banchoff's Morse Theories proven by Bloch \cite{Bloch13} this also explains the good correlation empirically observed \cite{RVRS} between the results using Forman's Morse Theory and Forman's Ricci curvature.

Clearly, the most important task ahead is to first translate these results and observations into implementable algorithms similar to those developed in \cite{HMR}, \cite{KSRS}, \cite{RVRS}, experiment with these on large scale real-life and model networks and hypernetworks, and compare the results with those in the works cited above. It is also necessary to compare between the methods proposed above, in particular if it is necessary to employ Bloch's more elaborate and less intuitive method, or if using Banchoff's approach, especially with the extensions to general cells and any direction for the height function suffices in practice.

On the theoretical plane, the immediate question arising is whether it is possible to extend Bloch's work from simplicial to more general cell complexes. 
Another problem would be to investigate whether, given the generalized versions of Banchoff's index, one can not dispense with the barycentric subdivision step in Bloch's proof of the equivalence of Forman's and Banchoff's Morse Theories, that is if it possible to apply Banchoff's curvature of general cells instead of that of vertices. (This would parallel Forman's 
assignment  of criticality to cells in any dimension, and not just to vertices.)
On a more general level, one would like to investigate if the polyhedral Morse Theory approaches render further results parallel to those holding in the classical theory (and partially for the Forman Morse Theory), such as the Morse inequalities and the Morse Theorem on degenerate critical points  (see \cite{Mil}).



\end{document}